\documentclass[a4paper]{article}

\usepackage[top=3cm, bottom=3cm, left=3cm, right=3cm]{geometry}

\usepackage[noadjust]{cite}
\usepackage{amsmath,amssymb,amsfonts,amsthm}
\usepackage{algorithmic}
\usepackage{graphicx}
\usepackage{textcomp}

\usepackage{stmaryrd}
\usepackage{pgf}
\usepackage{xcolor}
\usepackage{mathtools}
\usepackage{hyperref}
\usepackage{bbold}

\usepackage{array}

\usepackage{enumitem}

\newcommand{\R}{\mathbb{R}}
\newcommand{\diff}{\mathrm{d}}
\newtheorem{theorem}{Theorem}

\newtheorem{proposition}[theorem]{Proposition}

\theoremstyle{definition}
\newtheorem{definition}[theorem]{Definition}
\newtheorem{remark}[theorem]{Remark}

\renewcommand{\Re}{\operatorname{Re}}
\renewcommand{\Im}{\operatorname{Im}}
\DeclareMathOperator{\Var}{Var}
\DeclareMathOperator{\I}{I}
\DeclareMathOperator{\id}{id}

\begin{document}

\title{Strong stability of linear delay-difference equations}

\author{Felipe Gonçalves Netto\thanks{Universit\'e Paris Saclay, CNRS, CentraleSup\'elec, Inria, Laboratoire des Signaux et Syst\`emes, 91190 Gif-sur-Yvette, France.(E-mails: \url{felipe.goncalves-netto@centralesupelec.fr}, \url{guilherme.mazanti@centralesupelec.fr})}~\thanks{Department of Mathematics, University of Brasília, 70910-900 Brasília, Brazil.} \and Yacine Chitour\thanks{Universit\'e Paris Saclay, CNRS, CentraleSup\'elec, Laboratoire des Signaux et Syst\`emes, 91190 Gif-sur-Yvette, France. (E-mail: \url{yacine.chitour@centralesupelec.fr})}~\thanks{Fédération de Mathématiques de CentraleSupélec, 91190 Gif-sur-Yvette, France.} \and Guilherme Mazanti\footnotemark[1]~\footnotemark[4]}

\maketitle

\begin{abstract}
This paper considers linear delay-difference equations, that is, equations relating the state at a given time with its past values over a given bounded interval. After providing a well-posedness result and recalling Hale--Silkowski Criterion for strong stability in the case of equations with finitely many pointwise delays, we propose a generalization of the notion of strong stability to the more general class of linear delay-difference equations with an integral term defined by a matrix-valued measure. 
Our main result is an extension of Melvin Criterion for the strong stability of scalar equations, showing that local and global strong stability are equivalent, and that they can be characterized in terms of the total variation of the function defining the equation. We also provide numerical illustrations of our main result.

\bigskip

\end{abstract}

\section{Introduction}
\label{sec:introduction}

Time-delay systems are widely used to model phenomena that require considering a certain propagation time in their mathematical formulation for an accurate description. Such propagation time may come from several physical processes, such as the transportation of mass, energy, or information with finite speed, and may model communication lags, process duration, maturation times, reaction times, among other phenomena. In addition to their numerous applications, the mathematical analysis of time-delay systems raise nontrivial questions, which has motivated a large literature on this topic, in particular in Mathematics and Engineering \cite{Bellman1963Differential, Diekmann1995Delay, Hale1977Theory, Hale1993Introduction, Michiels2014Stability, Sipahi2011Stability}.

A particular class of time-delay systems are de\-lay-dif\-fe\-rence equations, sometimes also called con\-ti\-nuous-time difference equations or simply difference equations \cite[Chapter~9]{Hale1993Introduction}, \cite{Avellar1980Zeros, Chitour2016Stability, Cruz1970Stability, Hale1985Stability, Hale2002Strong, Henry1974Linear, Melvin1974Stability, Michiels2009Strong}. This class consists of those systems for which the value of the state at a given time $t$ can be expressed as a function of the state itself at previous times.

There are two main motivations for the analysis of delay-difference equations. On the one hand, properties of delay-difference equations, and in particular stability, are important in the analysis of the behavior of more general time-delay systems of neutral type, as detailed, for instance, in \cite[Chapter~9]{Hale1993Introduction}, \cite{Hale2002Strong, Henry1974Linear}, \cite[Section~1.4]{Mazanti2016Stability}. On the other hand, it has long been known that some classes of hyperbolic partial differential equations (PDEs) in one space dimension can be reduced, typically by using the method of characteristics, to delay-difference equations, and thus properties of the former can be studied through the analysis of the latter \cite{Chitour2016Stability}, \cite[Section~1.4.3]{Mazanti2016Stability}, \cite{BalogounNovel, Brayton1966Bifurcation, Chitour2021One, Chitour2017Persistently, Cooke1968Differential, Coron2008Dissipative, Fridman2010Bounds, Gugat2010Stars, Slemrod1971Nonexistence}. While in general such a transformation yields delay-difference equations with finitely many pointwise delays only (i.e., at time $t$, the state of the system depends on its values on finitely many points in the past depending on $t$), it has been established in \cite{Auriol2019Explicit} that some hyperbolic PDEs with in-domain coupling terms can be transformed into delay-difference equations with both pointwise and distributed delays (a distributed delay being a term involving an integral over past values of the state).

An important fact about the stability analysis of linear delay-difference equations with finitely many pointwise delays is that, even though exponential stability is preserved by perturbations on the matrices appearing in the description of the system, such a stability property can be lost under arbitrarily small perturbations of the delays \cite[Chapter~9, Section~6.1]{Hale1993Introduction}, \cite{Henry1974Linear, Melvin1974Stability, Moreno1973Zeros}. This has motivated the introduction of the notion of \emph{local strong stability}: a linear delay-difference equation with finitely many pointwise delays is said to be \emph{locally strongly stable} if it is exponentially stable and if it remains exponentially stable under small enough perturbations of the delays. Surprisingly, it turns out that this notion is equivalent to \emph{global strong stability} (sometimes named only \emph{strong stability}), which is defined as stability for \emph{any} positive values of delays. In addition, strong stability can also be characterized by a spectral condition depending only on the matrices appearing in the description of the system. These equivalences are known as the Hale--Silkowski Criterion for strong stability (which we recall in Theorem~\ref{HS:theorem} below), named after \cite{Hale1975Parametric, Silkowski1976Star}. We also note that the earlier work \cite{Melvin1974Stability} also provides some of those equivalences, but only for scalar equations, and this simpler case is known as Melvin Criterion.

In this paper, we are interested in delay-difference equations of the form
\begin{equation}\label{sys1}
x(t)=\int_{-1}^0 \diff M(\theta) x(t + \theta), 
\end{equation}
where $x(t)\in \R^n$ and $M\colon [-1,0]\to \mathcal{M}_n(\R)$ is a matrix-valued function with entries of bounded variation, assumed to be continuous from the right on $(-1,0)$, and the integral in the right-hand side of \eqref{sys1} is taken in the Riemann--Stieltjes sense. Note that the fact that we integrate over the interval $[-1, 0]$ in \eqref{sys1} is not a restriction, since the more general case where the integral in \eqref{sys1} is over an interval of the form $[-h, 0]$ for some $h > 0$ can be reduced to \eqref{sys1} by the change of time scale $y(t) = x(h t)$.

The delay-difference equation \eqref{sys1} can represent systems with both distributed and pointwise delays, by taking for instance $M$ to be the sum of a function in the Sobolev space $W^{1, 1}((0, 1), \mathcal M_n(\R))$ and a piecewise constant function. The main contributions of this paper are to extend notions of strong stability to systems of the form \eqref{sys1} and to provide a generalization of Melvin Criterion to \eqref{sys1} in the scalar case.

In \cite{Guang2011stability}, the notion of strong stability is considered for delay-difference equations \eqref{sys1} whose right-hand side is the sum of a finite number of pointwise delays and an integral term defined with $M\in W^{1, 1}((0, 1), \mathcal M_n(\R))$. However, the notion of strong stability used in \cite{Guang2011stability} only deals with perturbations in the pointwise delays, and not in the integral term. Instead, in this paper, the definition of strong stability we provide in Definition~\ref{def:strong-stab} also takes into account appropriate perturbations of the integral term.

The sequel of the paper is organized as follows. Section~\ref{sec:well-posed} deals with the well-posedness of \eqref{sys1}, providing an existence and uniqueness proof that strengthens a result of \cite{Hale1993Introduction} for the particular case of \eqref{sys1}. 
In Section~\ref{sec:strong}, after recalling the classical Hale--Silkowski Criterion for linear delay-difference equations with finitely many pointwise delays, we provide the definitions of strong stability for systems of the form \eqref{sys1}, and we state and prove our main result, Theorem~\ref{HS:scalar}, which is a generalization of Melvin Criterion to \eqref{sys1} in the scalar case. Numerical illustrations of our results are provided in Section~\ref{sec:numerical}, and a discussion closes the paper in Section~\ref{sec:concl}.

\medskip

\paragraph*{Notation.} In this paper, the set of $n \times n$ matrices with real coefficients is denoted by $\mathcal M_n(\mathbb R)$, and the identity matrix in $\mathcal M_n(\mathbb R)$ is denoted by $\I_n$, or simply $\I$ when $n$ is clear from the context. We use $\lvert \cdot \rvert$ to denote an arbitrary norm in $\mathbb R^n$ and $\lVert \cdot \rVert$ to denote its induced matrix norm in $\mathcal M_n(\mathbb R)$. The one-side limits of $f$ in $a$ are denoted by $f(a^+)$ for the right-sided limit and $f(a^-)$ for the left-sided limit. The indicator function of $A$ and the identity function in a given set are denoted by $\mathbb{1}_A$ and $\id$, respectively. The supremum norm of a continuous function defined on the interval $[-1,0]$ with values in $\R^n$ is denoted by $\Vert\cdot\Vert_\infty$.

\section{Well-posedness}
\label{sec:well-posed}

In this section we aim to provide global solutions for the delay-difference equation \eqref{sys1} in the space of continuous functions, according to the following definition.

\begin{definition}
A function $x$ is said to be a \emph{local solution} of \eqref{sys1} if there is $\delta > 0$ such that $x\in C([-1,\delta); \R^n)$ and $x$ satisfies \eqref{sys1} for $0\leq t < \delta$. We say that $x$ is a \emph{global solution} of \eqref{sys1} if it satisfies the definition of local solution with $\delta = +\infty$. In both cases, the restriction of $x$ to the interval $[-1, 0]$ is called the \emph{initial condition} of $x$.
\end{definition}

As usual with time-delay systems, given a function $x \in C([-1, \delta), \R^n)$ for some $\delta \in \R_+^\ast \cup \{+\infty\}$ and $t \in [0, \delta)$, the \emph{history function} of $x$ at time $t$ is the function $x_t \in C([-1, 0], \R^n)$ defined by
\[
x_t(\theta) = x(t+\theta), \quad \theta \in [-1, 0].
\]
With this notation, the initial condition of a solution $x$ of \eqref{sys1} is simply $x_0$. In the sequel, we set
\[
\mathcal C = C([-1, 0]; \R^n).
\]

\begin{remark}
If $x$ is a solution of \eqref{sys1}, then the initial condition $x_0$ satisfies
\[
x_0(0) = \int_{-1}^0 \diff M(\theta) x_0(\theta).
\]
We will thus consider in the sequel only initial conditions in the set $\mathcal C_0$ defined by
\[
\mathcal C_0 = \left\{\eta \in \mathcal C : \eta(0) = \int_{-1}^0 \diff M(\theta) \eta(\theta)\right\}.
\]
Note that, if $x \in C([-1, \delta); \R^n)$ is a solution of \eqref{sys1}, then $x_t \in \mathcal C_0$ for every $t \in [0, \delta)$.
\end{remark}

We will consider \eqref{sys1} for functions $M\colon [-1, 0] \to \mathcal M_n(\R)$ of bounded variation, according to the following definition.

\begin{definition}
Given a matrix-valued function $M\colon [a,b] \to \mathcal{M}_n(\R)$, we define the \emph{total variation} of $M$ in $[a, b]$ as
\[\Var M_{\rvert[a,b]}=\sup_P \sum_{i=1}^k \Vert M(t_{i+1})-M(t_{i}) \Vert,\]
where the supremum is taken over all partitions $P=\{[t_i,t_{i+1}]: 1\leq i \leq k\}$ of $[a,b]$. We say that $M$ is of \emph{bounded variation} in $[a, b]$ if $\Var M_{\rvert[a, b]} < +\infty$
\end{definition}

\begin{remark}
Even though $\Var M_{\rvert[a,b]}$ depends on the choice of the induced matrix norm $\lVert \cdot \rVert$, the notion of bounded variation is independent of the choice of the norm. In addition, $M\colon [a,b]\to \mathcal{M}_n(\R)$ has entries of bounded variation if and only if $\Var M_{\rvert [a,b]}<+\infty$. 
\end{remark}

\begin{remark}
\label{remk:M-mu}
Recall that, if $M$ is of bounded variation and $x$ is continuous, the integral in \eqref{sys1} is well-defined in the Riemann--Stieltjes sense. One can also see the integral in \eqref{sys1} as a Lebesgue--Stieltjes integral, by using the fact that there is a one-to-one correspondence between \emph{normalized} bounded variation functions $M\colon [-1,0]\to \mathcal M_n(\R)$ (where \emph{normalized} means that $M(-1)=0$ and $M$ is continuous from the right on $(-1,0)$) and matrix-valued Borel measures $\mu_M$ of bounded variation, expressed by
\[M(t)=\mu_M([-1,t]), \quad t \in [-1, 0].\]
In addition, denoting by $\lvert \mu_M \rvert$ the (scalar) total variation measure associated with $\mu_M$ and with respect to the induced norm $\lVert \cdot \rVert$ in $\R^n$, we have that $\Var M_{\rvert[-1, 0]} = \int_{-1}^0 \diff \lvert \mu_M \rvert$.

Note that the normalization $M(-1) = 0$ is not a restriction, since the integral in \eqref{sys1} does not change if a constant is added to $M$. From now on, we will use both Riemann--Stieltjes and Lebesgue--Stieltjes points of view, with the convention that matrix-valued functions of bounded variation will be denoted by capital Latin letters, while matrix-valued Borel measures of bounded variation will be denoted by lowercase Greek letters.

All the facts stated in this remark can be found in \cite[Appendix~I]{Diekmann1995Delay}.
\end{remark}

In order to state our result on existence and uniqueness of solutions of \eqref{sys1}, we define
\[
\mathcal W = \{M \colon [-1, 0] \to \mathcal M_n(\mathbb R) : M \text{ is of bounded variation and } \det(I - A_M) \neq 0\},
\]
where $A_M$ denotes the matrix $A_M = M(0) - M(0^-)$. According to Remark~\ref{remk:M-mu}, we also make the slight abuse of notation of saying that a matrix-valued Borel measure $\mu$ of bounded variation on $[-1, 0]$ belongs to $\mathcal W$ if the corresponding normalized bounded variation function belongs to $\mathcal W$.

We can now turn to our result on existence and uniqueness of solutions of \eqref{sys1}. The proof we present below is based on the existence proof from \cite[Chapter~2, Theorem~8.1]{Hale1993Introduction} for a more general neutral equation. Since our equation is linear, we can adapt that strategy to prove not only local existence but also uniqueness and global existence.

\begin{proposition}
Assume that $M \in \mathcal W$. Then, for every $\phi\in \mathcal C_0$, there exists a unique global solution $x\colon [-1,+\infty)\to \R^n$ of \eqref{sys1} with initial condition $x_0 = \phi$.
\end{proposition}

\begin{proof}
With no loss of generality, we assume that we have $\lim_{s\to 0^-} \Var M_{\rvert[s,0]} = 0$, for otherwise we can define
\[N(\theta)=\begin{cases}
(\I - A)^{-1} M(\theta), & \text{if } -1\leq \theta <0,\\
(\I - A)^{-1} M(0^-), & \text{if } \theta=0,
\end{cases}\]
with $A = M(0) - M(0^-)$, and remark that \eqref{sys1} is equivalent to
\[x(t)=\int_{-1}^0 \diff N(\theta)x(t+\theta),
\]
and $N$ satisfies $\Var N_{\rvert[s,0]}\to 0 $ as $s\to 0^-$.

Given $\alpha \in (0, 1]$, we define the set $\mathcal A(\alpha)$ by
\[
\mathcal A(\alpha) = \left\{\eta \in C([-1, \alpha]; \R^n) : \eta_0 = \phi\right\}
\]
and we endow $\mathcal A(\alpha)$ with the metric induced by the norm
\[
\lVert \eta \rVert_\alpha = \max_{s \in [-1, \alpha]} \lvert \eta(s) \rvert \quad \text{ in } C([-1, \alpha]; \R^n).
\]
We define the operator $S\colon \mathcal{A}(\alpha)\to\mathcal{A}(\alpha)$ by setting
\begin{equation}
\label{op:S}
S(x)(t) = \begin{dcases*}
\phi(t) & if $-1 \leq t \leq 0$, \\
\int_{-1}^0 \diff M(\theta) x(t + \theta) & if $0 < t \leq \alpha$,
\end{dcases*}
\end{equation}
and we remark that $S(x) \in \mathcal A(\alpha)$ for every $x \in \mathcal A(\alpha)$ since $\phi \in \mathcal C_0$. Notice also that $x$ is a local solution of \eqref{sys1} with initial condition $\phi$ if and only if $S(x) = x$.

For $x, y \in \mathcal{A}(\alpha)$, we have $S(x)(t) - S(y)(t) = 0$ for $-1 \leq t \leq 0$ and, for $0 < t \leq \alpha$, we have
\begin{align*}
    |S(x)(t) - S(y)(t)| & = \left|\int_{-1}^0 \diff M(\theta) \bigl[x(t + \theta) - y(t + \theta)\bigr] \right| \\
    & = \left| \int_{-t}^0 \diff M(\theta) \bigl[x(t + \theta) - y(t + \theta)\bigr]\right| \\
    & \leq \Var M_{\rvert[-t, 0]} \Vert x - y\Vert_{\alpha}\\
    & \leq \Var M_{\rvert[-\alpha, 0]} \Vert x - y\Vert_{\alpha}.
\end{align*}
Therefore, for every $x, y$ in $\mathcal{A}(\alpha)$, one has
\[\Vert S(x)-S(y)\Vert_{\alpha} \leq \Var M_{\rvert[-\alpha,0]} \Vert x - y \Vert_{\alpha}. \]
One deduces that, for some $\alpha > 0$ small enough that depends only on $M$, $S$ is a contraction. It therefore has a unique fixed point $x \in \mathcal{A}(\alpha)$, which is the unique local solution of \eqref{sys1} in $[-1, \alpha]$. As $\alpha > 0$ depends only on $M$ and not on $\phi$, one can iterate the above argument to extend the solution $x$ to $[-1, +\infty)$, obtaining a unique global solution.
\end{proof}

\section{Strong stability}
\label{sec:strong}

We now turn to the question of the strong stability of \eqref{sys1}. We start by recalling the definition of exponential stability.

\begin{definition}\label{def:ES}
We say that \eqref{sys1} is \emph{exponentially stable} if there are constants $\alpha, K > 0$ such that, for every solution $x$ of \eqref{sys1},
\[\Vert x_t\Vert_\infty \leq K e^{-\alpha t} \Vert x_0\Vert_\infty, \quad \text{for all } t\geq 0.\]
\end{definition}

We next recall the Hale--Silkowski Criterion in Section~\ref{sec:HS}, before providing our proposed extension of the strong stability notion for \eqref{sys1} in Section~\ref{sec:def-strong-stab} and finally stating and proving our main result, on the generalization of Melvin Criterion, in Section~\ref{sec:main-result}.

\subsection{Hale--Silkowski Criterion}
\label{sec:HS}

Let us consider the delay-difference equation
\begin{equation}\label{sys}
\Sigma(\tau)\colon \qquad x(t)=\sum_{k=1}^{N} A_k x(t-\tau_k),
\end{equation}
where $A =(A_1,\allowbreak \dotsc,\allowbreak A_N) \in \mathcal M_n(\R)^N$ and $\tau = (\tau_1,\allowbreak \dotsc,\allowbreak \tau_N) \in (0,1]^N$. Notice that \eqref{sys} corresponds to the particular case of \eqref{sys1} in which $M = \sum_{k=1}^N A_k \mathbb{1}_{[-\tau_k, 0]}$. We recall the following strong stability criterion for \eqref{sys}.

\begin{theorem}[Hale-Silkowski Criterion, \cite{Avellar1980Zeros, Hale1993Introduction, Silkowski1976Star}]\label{HS:theorem} Let $\bar\tau \in (0, 1]^N$ with rationally independent components and $\tau \in (0, 1]^N$. The following statements for \eqref{sys} are equivalent.

\begin{enumerate}[ref={\arabic*)}]
\item\label{item:HS-rat-indep} $\Sigma(\bar\tau)$ is exponentially stable.

\item\label{item:HS-rho} We have
\begin{equation}\label{HSC}
\rho_{\mathrm{HS}}(A) \overset{\mathrm{def}}{=} \max_{(\theta_1,\dots,\theta_N) \in [0, 2\pi]^N} \rho \left( \sum_{k=1}^N A_k e^{i\theta_k} \right) < 1,
\end{equation}
where $\rho (\cdot)$ denotes the spectral radius.

\item\label{item:HS-local} \emph{(Local strong stability)} There exists $\epsilon > 0$ such that, for every $r \in (0, 1]^N$  such that $|\tau - r| < \epsilon$, $\Sigma(r)$ is exponentially stable.

\item\label{item:HS-global} \emph{(Strong stability)} For every $r \in (0, 1]^N$, $\Sigma(r)$ is exponentially stable.
\end{enumerate}
\end{theorem}

Item \ref{item:HS-local} corresponds to the definition of local strong stability of \eqref{sys} at $\tau$, while \ref{item:HS-global} is the definition of global strong stability (or simply strong stability) of \eqref{sys}. The equivalence between \ref{item:HS-rho}, \ref{item:HS-local}, and \ref{item:HS-global} in Theorem~\ref{HS:theorem} for the scalar case (i.e., when $n = 1$) was provided first in \cite{Melvin1974Stability} and is known as Melvin Criterion, where the condition \eqref{HSC} becomes $\sum_{k=1}^N |A_k| < 1$.

\subsection{Strong stability for \texorpdfstring{\eqref{sys1}}{(\ref{sys1})}}
\label{sec:def-strong-stab}

Let us now provide our generalization of the notion of strong stability for \eqref{sys1}. We consider that a ``perturbation in the delays'' in \eqref{sys1} corresponds to a system of the form
\begin{equation}
\label{sys-perturbed}
x(t) = \int_{-1}^0 \diff M(\theta) x(t + \varphi(\theta)),
\end{equation}
where the function $\varphi\colon [-1, 0] \to [-1, 0]$ is a perturbation of the identity function. Note that, when $M = \sum_{k=1}^N A_k \mathbb{1}_{(-\tau_k, 0]}$ as in Section~\ref{sec:HS}, \eqref{sys-perturbed} becomes
\[
x(t) = \sum_{k=1}^N A_k x(t + \varphi(-\tau_k)),
\]
which is precisely a perturbation of the delays in \eqref{sys}.

If $\varphi$ is not continuous, the integral in \eqref{sys-perturbed} may fail to be defined in the Riemann--Stieltjes sense. However, denoting by $\mu$ the measure associated with $M$ through Remark~\ref{remk:M-mu}, under the assumption that $\varphi$ is Borel-measurable, \eqref{sys-perturbed} becomes
\begin{equation}
\label{sys2}
x(t) = \int_{-1}^0 \diff (\varphi_\ast\mu) (\theta)x(t+\theta),
\end{equation}
where $\varphi_\ast\mu$ is the pushforward of $\mu$ by $\varphi$ and the integral is in the Lebesgue--Stieltjes sense.

We are now in position to provide our generalization of the notion of strong stability.

\begin{definition}[Strong stability]
\label{def:strong-stab}
Let $M \in \mathcal W$, $\mu$ be the matrix-valued Borel measure associated with $M$ in the sense of Remark~\ref{remk:M-mu}, and define $\mathbf{B} = \{\varphi\colon[-1,0]\to [-1,0]\allowbreak : \allowbreak\varphi \text{ is}\allowbreak\text{Borel-measurable}\}$.
\begin{enumerate}[ref={\arabic*)}]
\item We say that \eqref{sys1} is \emph{locally strongly stable} if there exists $\epsilon>0$ such that, for all $\varphi \in \mathbf B$ with $\lVert \varphi - \id \rVert_\infty < \epsilon$  and $\varphi_\ast \mu \in \mathcal W$, \eqref{sys2} is exponentially stable.
\item We say that \eqref{sys1} is \emph{strongly stable} if, for all $\varphi\in\mathbf{B}$ such that $\varphi_\ast \mu \in \mathcal W$, \eqref{sys2} is exponentially stable.
\end{enumerate}
\end{definition}

\begin{remark}
It is not hard to verify that the condition $\varphi_\ast\mu \in \mathcal W$ is satisfied, in particular, if $M \in \mathcal W$ and $\varphi(\theta) \neq 0$ for all $\theta \in [-1, 0)$.
\end{remark}

\subsection{Strong stability criterion for \texorpdfstring{\eqref{sys1}}{(\ref{sys1})} in the scalar case}
\label{sec:main-result}

We are now in position to state and prove our main result, which generalizes Melvin Criterion to \eqref{sys1} in the scalar case.

\begin{theorem}\label{HS:scalar} Assume that $n = 1$. Then the following statements about \eqref{sys1} are equivalent.
\begin{enumerate}[ref={\arabic*)}]
\item\label{item:main-HS} The total variation of $\mu$ in $[-1, 0]$ is less than $1$, i.e.,
\begin{equation}\label{HSC_scalar}
\int_{-1}^0 \diff |\mu|(\xi)<1. 
\end{equation}
\item\label{item:main-loc} System \eqref{sys1} is locally strongly stable.
\item\label{item:main-glob} System \eqref{sys1} is strongly stable.
\end{enumerate}
\end{theorem}

\begin{proof}
It is clear that \ref{item:main-glob} implies \ref{item:main-loc}.

To show that \ref{item:main-loc} implies \ref{item:main-HS}, assume that
\[\int_{-1}^0 \diff |\mu|(\xi)\ge 1.\]
Recall that, by the Hahn and Jordan decomposition theorems (see, e.g., \cite[Theorems~3.3 and 3.4]{Folland1999Real}), there exist Borel-measurable subsets $P$ and $N$ of $[-1, 0]$ and unique positive measures $\mu^+$ and $\mu^-$ such that $\mu = \mu^+ - \mu^-$, $\lvert \mu\rvert = \mu^+ + \mu^-$, $\mu^-(P) = 0$, $\mu^+(N) = 0$, $P \cup N = [-1, 0]$, and $P \cap N = \emptyset$.

We need to show that, for every $\epsilon>0$, there exists 
$\varphi \in \mathbf B$ such that $\Vert\varphi - \id\Vert_{\infty}<\epsilon$, $\varphi_\ast\mu \in \mathcal W$, and \eqref{sys2} is not exponentially stable. So, given $\epsilon>0$, we take $Q=\{[t_k,t_{k+1}] : k \in \{0,\dotsc,m\},\, t_0=-1, t_{m+1}=0\}$ a partition of $[-1,0]$ such that
\[\lvert Q\rvert = \max_{k \in \{0, \dotsc, m\}} |t_{k+1}-t_k| < \epsilon.\]
Let us define the intervals $I_k=[t_k,t_{k+1})$ for $k\in \{0,\dotsc,m-1\}$ and $I_m=[t_m,t_{m+1}]$, and the sets
\[P_k=P\cap I_k\quad\text{and}\quad N_k=N\cap I_k.\]
We define the function $\varphi\colon [-1,0]\to [-1,0]$ by setting
\[
\varphi(\theta)=\begin{cases}
-\tau_k^P, & \text{if } \theta \in P_k \text{ for some } k,\\
-\tau_k^N, & \text{if } \theta \in N_k \text{ for some } k,
\end{cases}
\]
where, for $k \in \{0, \dotsc, m\}$, $-\tau_k^P$ and $-\tau_k^N$ are taken in $I_k$ and $\tau_0^P,\dotsc,\tau_m^P,\tau_0^N,\dotsc,\tau_m^N$ are rationally independent. By construction, we have $\varphi \in \mathbf B$,  $\varphi_\ast\mu \in \mathcal W$, $\lVert \varphi - \id \rVert_\infty < \epsilon$, and \eqref{sys2} becomes
\[
x(t) = \sum_{k=0}^m \left(a_k x(t - \tau_k^P) + b_k x(t - \tau_k^N)\right),
\]
where we have $a_k = \mu(P_k) = \mu^+(P_k) \geq 0$ and $b_k = \mu(N_k) = -\mu^-(N_k) \leq 0$.

Since $\sum_{k=0}^m \left(\lvert a_k\rvert + \lvert b_k\rvert\right) = \sum_{k=0}^m \left(\mu^+(P_k) + \mu^-(N_k)\right) = \lvert \mu \rvert ([-1, 0]) \geq 1$, we deduce that system \eqref{sys2} does not satisfy \eqref{HSC}. By Theorem~\ref{HS:theorem}, \eqref{sys2} is not exponentially stable, yielding the conclusion.

Assume now that \ref{item:main-HS} is satisfied and let us show \ref{item:main-glob}. Let $L = \int_{-1}^0 \diff \lvert\mu\rvert(\xi)<1$ and take $\varphi \in \mathbf B$. Any solution $x \in C([-1, +\infty); \R)$ of \eqref{sys2} satisfies
\begin{equation*}
|x(t)| = \left|\int_{-1}^0 x(t+\varphi(\theta))\, \diff\mu (\theta)\right| \leq \int_{-1}^0 |x(t+\varphi(\theta))|\, \diff |\mu|(\theta) 
\leq L \Vert x_t\Vert_\infty.
\end{equation*}
Hence, from \cite[Lemma~1]{Mazenc2015Trajectory}, \eqref{sys2} is exponentially stable, yielding that \ref{item:main-glob} holds true.
\end{proof}

\section{Numerical illustrations}
\label{sec:numerical}

To illustrate Theorem~\ref{HS:scalar}, we consider in this section System~\eqref{sys1} in the case where the measure $\mu$ associated with $M$ through Remark~\ref{remk:M-mu} is absolutely continuous with respect to the Lebesgue measure with an affine density, given by $\theta \mapsto a + b \theta$ for some real constants $a$ and $b$, in which case $M(\theta) = a(\theta + 1) + \frac{b}{2}(\theta^2 - 1)$ and \eqref{sys1} becomes
\begin{equation}\label{sys:a_b}
x(t) = \int_{-1}^0 x(t+\theta) (a+b\theta) \diff \theta.
\end{equation}
Our aim is to explicitly characterize the region of strong stability described by \eqref{HSC_scalar} in the $(a, b)$-plane, and also to numerically compute the region of exponential stability, in order to compare both.

\subsection{Region of strong stability}
\label{sec:example-strong}

Letting $\theta_v = -a / b$ to be the vertex of the parabola described by the function $M$ defined above (with the convention $\theta_v = 0$ if $b = 0$), we compute that $\Var M_{\rvert[-1,0]} = \lvert M(-1) - M(\theta_v) \rvert + \lvert M(\theta_v) - M(0)\rvert $ if $\theta_v \in (-1,0)$ and $\Var M_{\rvert[-1,0]} = \lvert M(-1) - M(0) \rvert$ otherwise. Noticing that $\theta_v \in (-1, 0)$ if and only if $a (b - a) > 0$, we thus obtain, by Theorem~\ref{HS:scalar}, the following result.

\begin{proposition}
System \eqref{sys:a_b} is strongly stable if and only if $(a, b)$ belongs to the union of the sets
\begin{align*}
\Bigl\{(a, b) \in \R^2 & : (a - b)^2 + a^2 < 2 \lvert b\rvert \text{ and }  a (b - a) > 0\Bigr\}, \\
\Bigl\{(a, b) \in \R^2 & : \bigl\lvert a - \frac{b}{2}\bigr\rvert < 1 \text{ and }  a (b - a) \leq 0\Bigr\}.
\end{align*}
\end{proposition}

\subsection{Region of exponential stability}
\label{sec:example-expo}

Let us now describe how we numerically approximate the region in the $(a, b)$-plane where \eqref{sys:a_b} is exponentially stable. From \cite[Chapter~12, Corollary~3.1]{Hale1977Theory}, we have that \eqref{sys:a_b} is exponentially stable if and only if there exists $\delta > 0$ such that all the zeros $\lambda$ of the characteristic equation
\[\Delta_{a,b}(s) = 1 - \int_{-1}^0 (a + b \theta) e^{s \theta} \diff \theta\]
satisfy $\Re \lambda \leq -\delta$. After a straightforward computation, we have that $\Delta_{a,b}(s) = \frac{Q(s)}{s^2}$, where
\[Q(s) = s^2 - a s + b + \big((a - b)s - b \big)e^{-s}.\]
Since $\Delta_{a, b}$ is holomorphic over $\mathbb C$, $s=0$ is a double root of $Q$, and the other roots of $Q$ coincide with the roots of $\Delta_{a,b}$.

We first notice that one can exclude a half-space from the exponential stability region of \eqref{sys:a_b}.

\begin{proposition}
System~\eqref{sys:a_b} is not exponentially stable on the half-space $\{(a, b) \in \R^2 : 2 - 2 a + b \leq 0\}$.
\end{proposition}

\begin{proof} 
We have $Q(0) = Q'(0) = 0$ and $Q''(0) = 2 - 2 a + b$, and thus $s=0$ is a triple root of $Q$ if $2-2a+b=0$, that is, $s=0$ is a root of $\Delta_{a,b}$ and \eqref{sys:a_b} is not exponentially stable. Now, if $Q''(0) = 2 - 2 a + b < 0$, we have that the restriction of $Q$ to the real line has a local maximum at $s=0$. So there is $s_1>0$ such that $Q(s_1) < 0$. On the other hand, 
$\lim_{s\to +\infty} Q(s)=+\infty$, so by continuity there is $s_2>s_1>0$ such that $Q(s_2)=0$. Hence $s_2$ is a positive real root of $\Delta_{a,b}$ and \eqref{sys:a_b} is not exponentially stable.
\end{proof}

We also remark that $Q$ can be seen as the characteristic function of a second-order time-delay system of retarded type and, as such, it admits only finitely many roots in any vertical strip of the complex plane (see, e.g., \cite[Corollary~1.9]{Michiels2014Stability}). In particular, \eqref{sys:a_b} is exponentially stable if and only if $\Re \lambda < 0$ for all roots $\lambda$ of $\Delta_{a, b}$. The function $Q$ admits the representation
\[
Q(s) = \det\left[s \I - A_0 - A_1 e^{-s}\right]
\]
where
\[A_0=\begin{pmatrix}
0 & 1\\ b & -a
\end{pmatrix} \text{ and } A_1=\begin{pmatrix}
0 & 0\\ -b & a-b
\end{pmatrix},\]
and thus, by \cite[Proposition~1.10]{Michiels2014Stability}, we have that, if $\lambda$ is a root of $Q$ with $\Re \lambda \geq 0$, then
\[|\lambda|\leq\alpha \overset{\mathrm{def}}{=} \Vert A_0\Vert_2 + \Vert A_1\Vert_2.\]
Taking into account that roots of $Q$ appear in complex conjugate pairs, we deduce that, given any root $\lambda$ of $Q$ with nonnegative real part, either $\lambda$ or $\bar\lambda$ belong necessarily to the square $S_\alpha = \{\lambda \in \mathbb C : 0 \leq \Re \lambda \leq \alpha,\, 0 \leq \Im \lambda \leq \alpha\}$. 

In order to numerically approximate the region of exponential stability of \eqref{sys:a_b} in the $(a, b)$-plane, we thus discretize the region of interest into a grid $(a_j, b_k)$, for $j$ and $k$ in suitable finite sets. For each pair $(a_j, b_k)$, the system is not exponentially stable if $2 - 2 a_j + b_k \leq 0$. Otherwise, we compute numerically the roots of $Q$ in the square $S_\alpha$ using the QPmR algorithm from \cite{Vyhlidal2014QPmR}\footnote{More precisely, we use the Python implementation available at \url{https://github.com/DSevenT/QPmR}.}. If $Q$ admits roots in this square with nonnegative real part other than the double root at $0$, then the system is not exponentially stable, otherwise it is exponentially stable.

\subsection{Comparisons}

\begin{figure*}[ht]
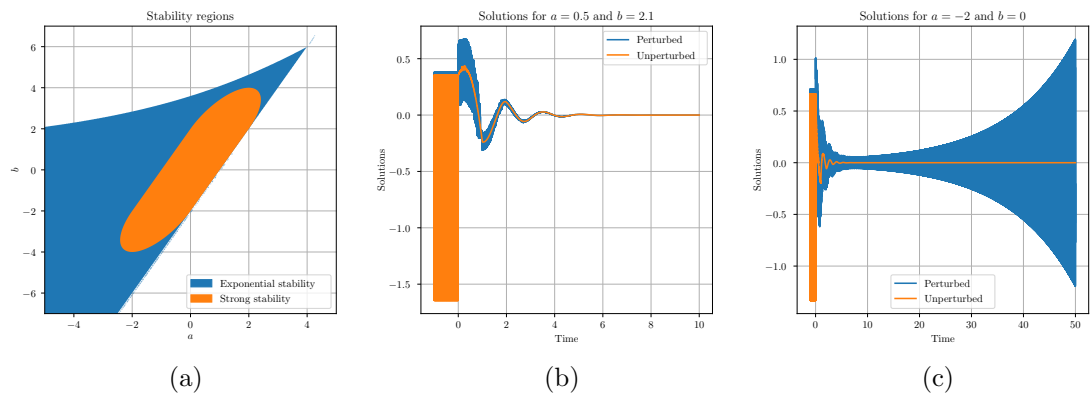

\centering
\begin{tabular}{@{} >{\centering} m{0.33\textwidth} @{} >{\centering} m{0.33\textwidth} @{} >{\centering} m{0.33\textwidth} @{}}
\resizebox{0.33\textwidth}{!}{\input{Figures/stability.pgf}} & \resizebox{0.33\textwidth}{!}{\input{Figures/trajectories_strongly_stable.pgf}} & \resizebox{0.33\textwidth}{!}{\input{Figures/trajectories_exponentially_stable.pgf}} \tabularnewline
(a) & (b) & (c)
\end{tabular}
\caption{(a) Stability regions and (b, c) solutions of \eqref{sys:a_b} (orange curves) and of \eqref{sys:a_b_perturbed} (blue curves) (b) for the pair $(a, b) = (0.5, 2.1)$ in the strong stability region and (c) for the pair $(a, b) = (-2, 0)$ in the exponential stability region, but outside of the strong stability region.}
\label{stab_regions}
\hrulefill
\end{figure*}

Fig.~\ref{stab_regions}(a) represents the regions of exponential stability (in blue) and of strong stability (in orange), computed according to the discussions in Sections~\ref{sec:example-strong} and \ref{sec:example-expo}.

In order to illustrate the effects of the lack of strong stability, Fig.~\ref{stab_regions}(b) and (c) provide numerical simulations of solutions of \eqref{sys:a_b} and of the perturbed system
\begin{equation}
\label{sys:a_b_perturbed}
x(t) = \int_{-1}^0 x(t + \varphi(\theta)) (a + b \theta) \diff\theta.
\end{equation}
We selected a perturbation $\varphi\colon [-1, 0] \to [-1, 0]$ which is constant in each interval of the form $(-\frac{j}{N}, -\frac{j-1}{N})$ for $j \in \{1, \dotsc, N\}$ and with $N = 25$, with randomly chosen values in each of these intervals in such a way that $\lVert \varphi - \id\rVert_\infty \leq \frac{1}{N} = 0.04$. As the lack of strong stability usually results from the apparition of unstable roots with large imaginary part in the perturbed system, we selected an initial condition of the form $x_0(t) = \cos(2 k \pi t) + c_0$, with $k = 127$ and $c_0 \in \mathbb R$ chosen so that $x_0 \in \mathcal C_0$. We observe, in Figure~\ref{stab_regions}(b) and (c), that the corresponding solutions of \eqref{sys:a_b} converge exponentially to the origin in both cases $(a, b) = (0.5, 2.1)$ and $(a, b) = (-2, 0)$, as expected, but, for the perturbed system \eqref{sys:a_b_perturbed}, the corresponding solution converges to the origin in the strongly stable case $(a, b) = (0.5, 2.1)$, but it diverges in the case $(a, b) = (-2, 0)$ where strong stability does not hold true.

\section{Conclusions}
\label{sec:concl}

In this paper, we investigate strong stability for a linear system of delay-difference equations. We extend the concept of strong stability to a more general class of systems which may involve in particular distributed delays. We present the well-posedness of the $n$-dimensional case and provide a criterion for strong stability in the scalar case, which generalizes Melvin Criterion. Additionally, we present numerical illustrations for a particular system in order to visualize and compare the regions of exponential stability and strong stability. Future work will consist on deriving strong stability criteria for \eqref{sys1} in the multi-dimensional case.

\bibliographystyle{abbrv}
\bibliography{root}

\end{document}